\documentclass[12pt,a4paper]{article}
\usepackage[utf8]{inputenc}
\usepackage[english]{babel}
\usepackage[T1]{fontenc}
\usepackage{amsfonts, amssymb, mathrsfs, amsmath, amsthm}
\usepackage{graphicx,psfrag,epsf}
\usepackage{fixltx2e}
\usepackage{xcolor}
\usepackage{natbib}
\usepackage{geometry}
\usepackage{calc}
\usepackage{tikz}
\usepackage{subfig}
\usepackage{enumerate}
\usepackage{url} 
\usepackage{soul} 
\usepackage{array}
\usepackage{cleveref} 

\title{A new inequality for maximum likelihood estimation in statistical models with latent variables}
\author{Niels Lundtorp Olsen \\ Department of Applied Mathematics and Computer Science, \\ Technical Unversity of Denmark}

\newcolumntype{L}{>{\centering\arraybackslash}m{1.5cm}}

\newcommand{\de}{\: \mathrm{d}}

\theoremstyle{definition}

\newtheorem{thm}{Theorem}

\newtheorem*{defi*}{Definition}

\begin{document}
\maketitle

\begin{abstract}
	Maximum-likelihood estimation (MLE) is arguably the most important tool for statisticians, and many methods have been developed to find the MLE.
	We present a new inequality involving posterior distributions of a latent variable that holds under very general conditions. It is related to the EM algorithm and has a clear potential for being used in a similar fashion.
\end{abstract}

\section{Introduction} Since and before R.A. Fisher introduced the concept of likelihood and maximum-likelihood estimation \citep{fisher1922}, it has been a goal to find the maximum likelihood for a given statstical model, which in principle is done solving the associated \emph{score equations}.
Maximum-likelihood estimation is arguably the most used principle for statistical inference and is underpinned by a lot of theory. 

However, due to the fact that solving the score equations is often infeasible,
 auxillary methods have been developed specifically to facilitate maximization of the likelihood function, most notably the \emph{{Expectation}-{Maximi\-zation} \emph{(EM)} algorithm} \citep{dempster1977}. The huge increase in computational power during the last decades has also faciliated many new methods, in particular simulation-based approaches.

\paragraph{Latent variables}
A notable challenge in maximum likelihood estimation (and indeed, any statistical inference) is the presence of \emph{latent} or \emph{unobserved} variables. A latent variable $w$ is characterized by the fact that it acts as part of the statistical model, but is not observed. In terms of maximum likelihood inference, this implies the presence of a sum (if $w$ is discrete) or an integral (if $w$ is continuous) in the likelihood expression. Intergrals are notoriously difficult to evaluate, so other approaches are often needed. 

There are three common approaches to handle latent variables in maximum likelihood estimation (sometimes in combination):

\begin{itemize}
	\item EM algorithm and its derivatives
	\item Monte Carlo methods
	\item Approximation methods, such as \emph{Laplace's method}.
\end{itemize}
In this article, we present and proof a new inequality involving latent variables, that when true guarantees an increase in likelihood. The resulting theorem is related to the EM algorithm as it also uses the posterior of the latent variable, but our inequality is more general and does not imply an algorithm per se. 
We will briefly discuss some practicalities but otherwise leave applications for future work. 

\section{Theorem}

Suppose that we are given a statistical family consisting of an observation $y \in \mathcal{Y}$, a latent variable $w \in W$ and an unknown parameter $\theta$ in parameter space $\Theta$. 

Assume that the joint variable $(y,w)$ is dominated; that is $P_\theta((y, w) \in A) = \int_A p_{\theta}(y, w) d (\lambda \otimes \mu)(y, w)$ for a measure $\lambda$ on $\mathcal{Y}$ and $\mu$ on $\mathcal{W}$, and assume that for every $\theta$, $p_{\theta}(y, w)$ is non-zero for almost all $w, y$.

Let  $p_\theta(y)$ be the marginal density for $y$, and let $L(\theta) := p_\theta(y),L(\theta, w) := p_\theta(y, w) $ denote marginal and posterior likelihoods, respectively.
	
\begin{thm}
	Let $\theta_1, \theta_2 \in \Theta$. 
Then $L(\theta_1) < L(\theta_2)$ iff the following the following inequality is true:
\begin{align}
	\int \min\left(1, \frac{L(\theta_2 , w)}{L(\theta_1 , w)} \right) \de P_{\theta_1}(w | y) > 
	\int \min\left(1, \frac{L(\theta_1 , w)}{L(\theta_2 , w)} \right) \de P_{\theta_2}(w | y),
\end{align}	
where $P_{\theta}(w | y)$ is the posterior distribution of $w$ under $\theta$ given $y$. 
\end{thm}
That is, by integrating the "truncated likelihood-ratios" under the posterior distributions, we can compare the likelihood of $y$ under $\theta_1$ and $\theta_2$. 

\begin{proof}
	Let $A$ denote the subset of $W$ where $L(\theta_2 , w) <  L(\theta_1 , w)$. 
	First consider the left integral:
\begin{multline}
	\int \min\left(1, \frac{L(\theta_2 , w)}{L(\theta_1 , w)} \right) \de P_{\theta_1}(w | y) = 
	\int \min\left(1, \frac{L(\theta_2 , w)}{L(\theta_1 , w)} \right)  \frac{L(\theta_1, w)}{L(\theta_1)} \de \mu(w) = \\ \int_A \frac{L(\theta_2 , w)}{L(\theta_1)} \de \mu(w) + \int_{A^c} \frac{L(\theta_1 , w)}{L(\theta_1)} \de \mu(w) = \\
	\frac{1}{L(\theta_1)}  \int 1_{A}(w) L(\theta_2 , w) + 1_{A^c}(w) L(\theta_1 , w) \de \mu(w)
\end{multline}	
We get a similar result for the right integral with $L(\theta_1)$ replaced by $L(\theta_2)$. Now the theorem follows.
\end{proof}

\subsection*{Remarks}
Note that the truncated likelihood-ratio is numerically stable due to the upper limit of 1. 

\paragraph{Application of the theorem} In general, the posterior distributions $P_{\theta_1}(w | y)$ and  $P_{\theta_2}(w | y)$ are not attainable. Therefore, Monte Carlo methods would have to be applied.

Note that the theorem does not imply an algorithm per se. However, by identifying $\theta_1$ with the current estimate and $\theta_2$ with a proposed estimate, we can outline an algorithm, if we can come up with new proposed estimates. This  will obviously be dependent on the model in question, and we leave this for future work.
	
\section{Discussion}

With its very general setting, the presented theorem should be applicable in a wide range of models, since latent variables are present in many classes of statistical models.

An interesting notion is the fact the theorem does not require any assumptions on the parameter space $\Theta$. 
The most important assumption is that $p_\theta(y, w)$ is nonzero for almost all $w$, which can actually be relaxed to some extent.

\subsection*{Comparison to the EM algorithm}
A close relative is the {EM} algorithm and its derivatives -- using the current estimate $\theta_0$, the EM algorithm proposes a new estimate $\theta^*$, that is guaranteed to improve the likelihood. 

However, there are some notable differences:
\begin{itemize}
	\item The EM algorithm and its derivatives use only the posterior distribution of the current estimate.
	\item Application of the theorem requries using two posterior distributions. 
\end{itemize}
We believe the great potential of the presented theorem is when combined with a clever 
"proposal scheme" for new estimates. 
 Unlike the EM algorithm and its derivatives, we are free to choose any propsal $\theta^*$ for an updated estimate. Since the proposal scheme can be tailored to a specific model, it should be possible to create new and powerful methods to facilitate maximum likelihood estimation.

\bibliographystyle{agsm}

\bibliography{bib}

\end{document}